\documentclass[12pt]{article}
\usepackage{amsfonts,amsthm}
\usepackage{amsmath,lineno}
\usepackage{graphicx,epsfig}
\usepackage{amssymb}
\usepackage{subfigure,color,setspace}

\DeclareMathOperator{\sgn}{sgn}
\begin{document}

\title{Orhonormal wavelet bases on the 3D ball via volume preserving map from the regular octahedron}

\author{Adrian Holho\c s, Daniela Ro\c sca}
%
%
%
\maketitle
\begin{abstract}
We construct a new volume preserving map from the unit ball $\mathbb B^3$ to the regular 3D octahedron, both centered at the origin, and its inverse. This map will help us to construct refinable grids of the 3D ball, consisting in diameter bounded cells having the same volume. On this 3D uniform grid, we construct a multiresolution analysis and orthonormal wavelet bases of $L^2(\mathbb B^3)$, consisting in piecewise constant functions with small local support. 
\end{abstract}
Keywords: equal volume projection, hierarchical grid, wavelets\\
MSC 2010: 65T60, 65N55
\def \unu #1{#1^{\prime}}
\section{Introduction}
\label{intro}

Spherical 3D signals occur in a wide range of fields, including  computer graphics, medical imaging, cristallography (texture analysis of crystals) and geoscience (geography, geodesy, meteorology, seismology, astronomy, etc). Therefore, we need suitable efficient techniques for manipulating such signals, and one of the most efficient technique consists in using wavelets on the 3D ball (see e.g. \cite{Englezii, flag, Mich} and the references therein). In this paper we propose a construction orthonormal basis of wavelets with small support, defined on the 3D ball $\mathbb B^3$, starting from a multiresolution analysis. Our wavelets will be piecewise constant functions on the cells of a uniform and refinable grid of $\mathbb B^3$. By a refinable (or hierarchical) grid we mean that  the cells can be divided successively into a given number of smaller cells of the same volume. By a uniform grid we mean that all the cells at a certain level of subdivision have the same volume. These two very important properties of our grid derive from the fact that it is constructed by mapping a uniform and refinable grid of the 3D regular octahedron,  using a volume preserving map onto $\mathbb B^3$. 
Compared to the wavelets on the 3D ball constructed in \cite{flag} and \cite{Mich}, with localized support, our wavelets have local  support, and this is very important when dealing with data consisting in big jumps on small portions, as shown in \cite{rosmc}. Another construction of piecewise constant wavelets on the 3D ball was realized in \cite{Englezii}, starting from a similar construction on the 2D sphere. The author pretend that his wavelets are the first Haar wavelets on the 3D ball which are orthogonal and symmetric, even though we do not see any symmetry, neither in the cells, nor in the decomposition matrix. Moreover, his $8\times 8$ decomposition matrices change in each step of the refinement, the entries depending on the volumes of the cells, which are, in our opinion, difficult to evaluate and for this reason they are not calculated expicitely in \cite{Englezii}.  Another advantage of our construction is that our cells are diameter bounded, unlike the cells in \cite{Englezii} containing the origin, which become long and thin after some steps of refinement.

The paper is structured as follows. In Section \ref{section2} we introduce some notations used for the construction of the volume preserving map. In Section \ref{section3} we construct the volume preserving maps between the regular 3D octahedron and the  3D ball $\mathbb B^3$. In Section \ref{section4} we construct a uniform refinable grid of the regular octahedron followed by implementation issues, and its projection onto  $\mathbb B^3$ . Finally, in Section \ref{section5} we construct a multiresolution analysis and piecewise constant wavelet bases of $L^2( \mathbb B^3)$.
\section{Preliminaries}\label{section2}
Consider the ball of radius $r$ centered at the origin $O$, defined as
\[\mathbb B^3 =\left\{(x,y,z)\in\mathbb R^{3}, x^{2}+y^2+z^2\leq r^2\right\}\]and the regular octahedron $\mathbb K$ of the same volume, centered at $O$ and with vertices on the coordinate axes
\[\mathbb K=\left\{(x,y,z)\in\mathbb R^3, |x|+|y|+|z|\leq a\right\}.\] Since the volume of the regular octahedron is ${4a^{3}}/3$,  we have
\begin{equation}\label{a}
a=r\sqrt[3]\pi.
\end{equation}
The parametric equations of the ball are
\begin{eqnarray}
&&x=\rho \cos \theta \sin \varphi,\nonumber\\
&&y=\rho\sin \theta \sin \varphi,\label{parsf}\\
&&z=\rho \cos \varphi,\nonumber
\end{eqnarray}
where $\varphi \in [0,\pi]$ is the colatitude, $\theta \in [0,2\pi)$ is the longitude and $\rho\in[0,r]$ is the distance to the origin. A simple calculation shows that the volume element of the ball is
\begin{equation}\label{ve}
dV=\rho^2 \sin \varphi\, d\rho\,d\theta\,d\varphi.
\end{equation}
 
The ball and the octahedron can be split into eight congruent parts, each part being situated in one of the eight octants $I_{i}^\pm$, $i=0,1,2,3$,
\begin{eqnarray*}
I_0^+=\{ (x,y,z),\ x\geq 0,\ y\geq 0,\ z\geq 0\},\quad
I_0^-=\{ (x,y,z),\ x\geq 0,\ y\geq 0,\ z\leq 0\},\\
I_1^+=\{ (x,y,z),\ x\leq 0,\ y\geq 0,\ z\geq 0\},\quad
I_1^-=\{ (x,y,z),\ x\leq 0,\ y\geq 0,\ z\leq 0\},\\
I_2^+=\{ (x,y,z),\ x\leq 0,\ y\leq 0,\ z\geq 0\},\quad
I_2^-=\{ (x,y,z),\ x\leq 0,\ y\leq 0,\ z\leq 0\},\\
I_3^+=\{ (x,y,z),\ x\geq 0,\ y\leq 0,\ z\geq 0\},\quad I_3^-=\{
(x,y,z),\ x\geq 0,\ y\leq 0,\ z\leq 0\}.
\end{eqnarray*}

Let $\mathbb B_{i}^s$ and $\mathbb{K}_i^s$ be the regions of $\mathbb B^{3}$ and $\mathbb K$, situated in $I_{i}^s$, respectively.

\begin{figure}
\centering
\includegraphics[width=0.35\textwidth]{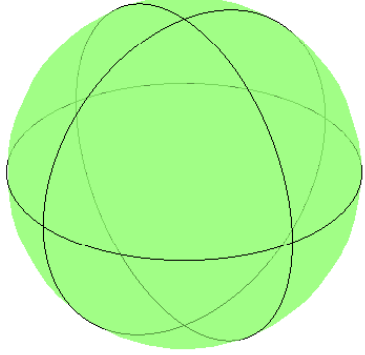}
\caption{The eight spherical zones obtained as intersections of the coordinate planes with the ball $\mathbb B^3$.}
\label{diag}
\end{figure}

Next we will construct a map $\mathcal U: \mathbb B^3\to \mathbb
K$ which preserves the volume, i.e. $\mathcal U$ satisfies
\begin{equation}\label{app}
\mbox{vol}(D)=\mbox{vol}(\mathcal U(D)),\qquad \mbox{for all }D\subseteq \mathbb B^3,
\end{equation}
where $\mbox{vol}(D)$ denotes the volume of a domain $D.$ For an arbitrary point $(x,y,z)\in \mathbb B^3$ we denote
\begin{equation}\label{notu}
(X,Y,Z)=\mathcal U(x,y,z) \in \mathbb K.
\end{equation}
\section{Construction of the volume preserving map $\mathcal U$ and its inverse}\label{section3}

\begin{figure}
\begin{center}
\includegraphics[width=0.4\textwidth]{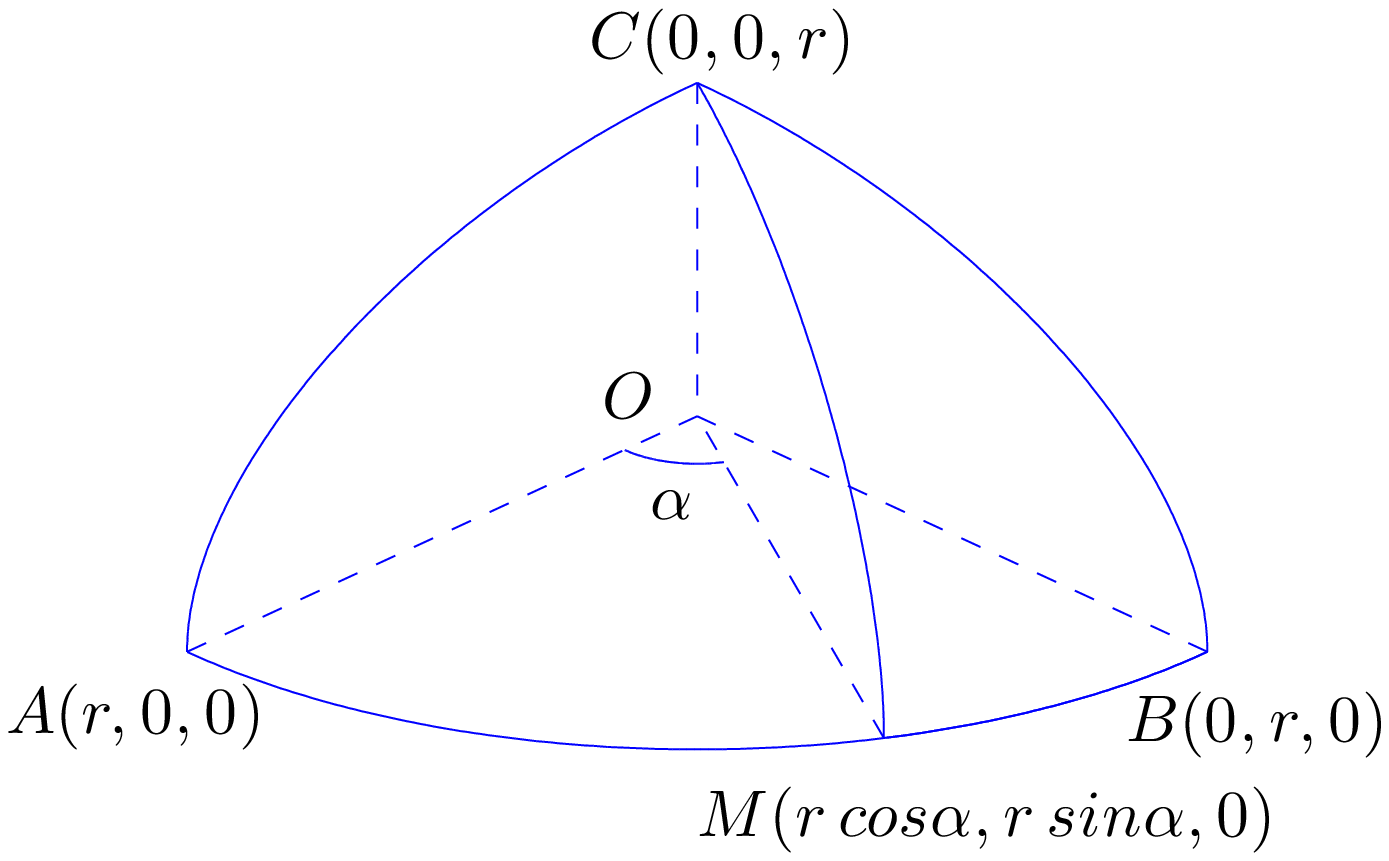}\hspace{1cm}
\includegraphics[width=0.4\textwidth]{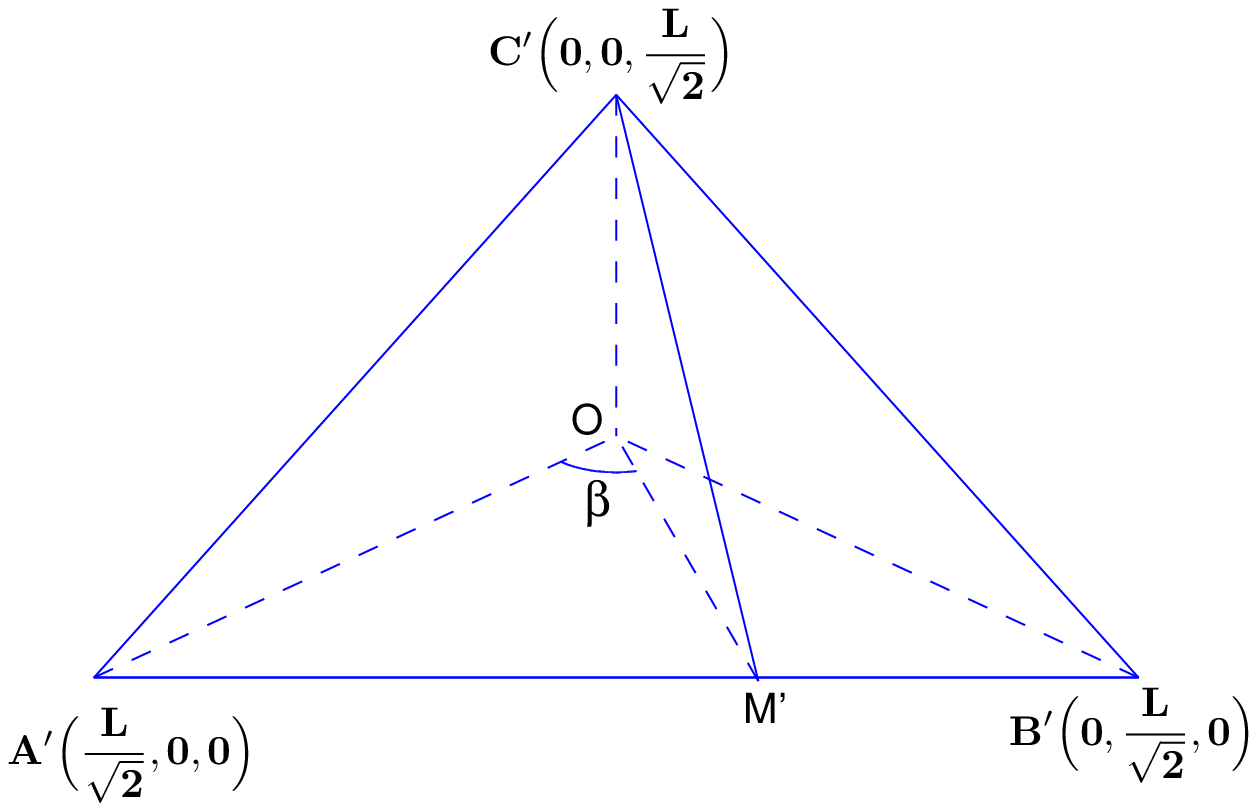}
\end{center}
\caption{The spherical region ${OAMC}$ and its image $OA'M'C'=\mathcal U(OAMC)$ on the octahedron. }\label{desexpl}
\end{figure}

We focus on the region $\mathbb B_0^+\subset I_{0}^+$ where we
consider the points $A=(r,0,0)$, $B=(0,r,0)$, $C=(0,0,r)$ and the
vertical plane of equation $y=x\tan \alpha$ with $\alpha \in
(0,\pi/2)$ (see Figure \ref{desexpl} (left)). We denote by $M$ its
intersection with the great arc $\widetilde{AB}$ of the sphere of
radius $r$. More precisely, $M=(r\cos \alpha, r\sin \alpha, 0).$
The volume of the spherical region $OAMC$ equals $r^3\alpha/3.$
Now we intersect the region $\mathbb K_1^+$ of the octahedron with
the vertical plane of equation $y=x \tan \beta$ and denote by
$M'(m,n,0)$ its intersection with the edge $A'B'$, where $A'\left(
a,0,0 \right),$ $B'\left( 0,a,0 \right)$ (see Figure \ref{desexpl}
(right)). Then $m+n=a$ and from $n=m \tan \beta$ we find
$$
m=a\cdot \frac 1{1+\tan \beta},\qquad n=a\cdot\frac {\tan \beta} {1+\tan \beta}.
$$

The volume of $OA'M'C'$ is
\[\mathcal{V}(OA'M'C')=\frac{OC'\cdot \mathcal{A}(OA'M')}{3}=\frac{a}{3}\cdot \frac{OA'\cdot n}{2}=\frac{a^3 \tan\beta}{6(1+\tan \beta)}.\]
If we want the volume of the region $OAMC$ of the unit ball to be equal to the volume of $OA'M'C'$, we obtain
$$
\alpha=\frac{\pi}2\cdot\frac{\tan \beta}{1+\tan \beta},\mbox{ whence }\tan \beta=\frac{2 \alpha}{\pi-2\alpha}.
$$
This give us a first relation between $(x,y,z)$ and $(X,Y,Z)$:
\begin{equation*}
\frac{Y}{X}=\frac{2 \arctan \frac{y}{x}}{\pi-2 \arctan \frac{y}{x}}.
\end{equation*}
Using the spherical coordinates (\ref{parsf}) we obtain
\begin{equation}\label{rel1}
Y= \frac{2\theta}{\pi-2\theta}\cdot X.
\end{equation}
In order to obtain a second relation between $(x,y,z)$ and $(X,Y,Z)$, we impose that, for an arbitrary $\widetilde \rho\in (0,r]$ the region
\[\left\{(x,y,z)\in\mathbb R^3,\ x^2+y^2+z^2\leq \widetilde \rho^2,\ x,y,z\geq 0 \right\}
\mbox{ of volume }\frac{\pi \widetilde \rho^3}6
\] is mapped by $\mathcal U$ onto
\[\left\{(X,Y,Z)\in\mathbb R^3,\ X+Y+Z\leq \ell,\ X,Y,Z\geq 0\right\}
\mbox{ of volume }\frac{\ell^3}6
.\] Then, the volume preserving condition \eqref{app} implies $\ell= a\cdot \widetilde \rho/r$, with $a$ satisfying \eqref{a}. Thus,
\[X+Y+Z=\frac a r \sqrt{x^2+y^2+z^2}\]
and in spherical coordinates this can be written as
\begin{equation}\label{rel2}
X+Y+Z=\frac{a\rho}r.
\end{equation}
In order to have a volume preserving map, the modulus of the
Jacobian $J(\mathcal U)$ of $\mathcal U$ must be 1, or,
equivalently, taking into account the volume element \eqref{ve},
we must have
\begin{equation} \label{rel3} J(\mathcal U)=\begin{vmatrix}\unu{X}_{\rho} & \unu{X}_{\varphi} & \unu{X}_{\theta} \\ \unu{Y}_{\rho} & \unu{Y}_{\varphi} & \unu{Y}_{\theta}\\ \unu{Z}_{\rho} & \unu{Z}_{\varphi} & \unu{Z}_{\theta} \end{vmatrix}=\rho^2\sin\varphi.
\end{equation}
Taking into account formula \eqref{rel2}, we have
\[J(\mathcal U)=\begin{vmatrix}\unu{X}_{\rho} & \unu{X}_{\varphi} & \unu{X}_{\theta} \\ \unu{Y}_{\rho} & \unu{Y}_{\varphi} & \unu{Y}_{\theta}\\ a/r-\unu{X}_{\rho}-\unu{Y}_{\rho} & -\unu{X}_{\varphi}-\unu{Y}_{\varphi} & -\unu{X}_\theta -\unu{Y}_\theta \end{vmatrix}=\begin{vmatrix}\unu{X}_{\rho} & \unu{X}_{\varphi} & \unu{X}_{\theta} \\ \unu{Y}_{\rho} & \unu{Y}_{\varphi} & \unu{Y}_{\theta}\\ a/r & 0 & 0 \end{vmatrix}=\frac a r \begin{vmatrix} \unu{X}_{\varphi} & \unu{X}_{\theta} \\  \unu{Y}_{\varphi} & \unu{Y}_{\theta} \end{vmatrix}. \]
Further, using relation (\ref{rel1}) we get
\[J(\mathcal U)=\frac ar \begin{vmatrix} \unu{X}_{\varphi} & \unu{X}_{\theta} \\  \frac{2\theta}{\pi-2\theta}\cdot\unu{X}_{\varphi} & \frac{2\theta}{\pi-2\theta}\cdot\unu{X}_{\theta}+ \frac{2\pi}{(\pi-2\theta)^2}\cdot X \end{vmatrix}=\frac ar \begin{vmatrix} \unu{X}_{\varphi} & \unu{X}_{\theta} \\  0 & \frac{2\pi }{(\pi-2\theta)^2}\cdot X \end{vmatrix}= \frac{2\pi a}{r(\pi-2\theta)^{2}}XX'_{\varphi}.\]
For the last equality, we have multiplied the first row by
$-2\theta/(\pi-2\theta)$ and we have added it to the second row.
Then, using the expression for $J(\mathcal U)$ obtained in
\eqref{rel3} we get the differential equation
\[2\unu{X}_{\varphi}\cdot X=\frac{r\rho^2}{\pi a  } (\pi-2\theta)^2\sin\varphi .\]

The integration with respect to $\varphi $ gives
$$
X^{2}=-\frac {r(\pi-2\theta)^{2}}{\pi a}\rho^{2}\cos \varphi+\mathcal C(\theta,\rho),
$$
and further, for finding $\mathcal C(\theta,\rho)$ we use the fact that, for $\varphi=\pi/2$ we must obtain $Z=0$. Thus, for $\varphi=\pi/2$ we have
$$
X^{2}=\mathcal C(\theta,\rho), \mbox{ so }Y=\frac{2\theta}{\pi-2\theta}\sqrt{\mathcal C(\theta,\rho)}, \mbox{ and }
$$
$$
Z=\frac{a\rho}r-X-Y=\frac{a\rho}r-\frac{\pi}{\pi-2\theta}\sqrt{\mathcal C(\theta,\rho)}.
$$
Thus, $Z=0$ is obtained for
$$
\mathcal C(\theta,\rho)=\frac{a^2\rho^2}{\pi^{2}r^2}(\pi-2\theta)^{2},
$$
and finally, the map $\mathcal{U}$ restricted to the region $I_{0}^+$ is

\begin{flalign}
X&=\frac{\sqrt 2 }{\pi^{2/3}}\cdot \rho (\pi-2\theta)\sin \frac \varphi 2,\label{scX}\\
Y&=\frac{\sqrt 2}{\pi^{2/3}}\cdot \rho \cdot 2\theta \sin \frac \varphi 2,\label{scY}\\
Z&=\pi^{1/3}\rho \big(1-\sqrt 2\sin \frac \varphi 2\big)\label{scZ}.
\end{flalign}

In the other seven octants, the map $\mathcal U$ can be obtained by symmetry as follows. A point $(x,y,z)\in \mathbb B^3,$ can be written as
$$
(x,y,z)=(\sgn(x)\cdot|x|, \sgn(y)\cdot|y|,\sgn(z)\cdot|z|),\quad \mbox{ with }(|x|,|y|,|z|)\in I_{0}^+.
$$
Therefore, if we denote by $(\overline X,\overline Y,\overline Z)=\mathcal U(|x|,|y|,|z|)$, then we can define $\mathcal U(x,y,z)$ as
\begin{equation}\label{prelsim}
\mathcal U(x,y,z)=(\sgn(x)\cdot \overline X,\sgn(y)\cdot \overline Y,\sgn(z)\cdot \overline Z).
\end{equation}

Next we deduce the formulas for the inverse of $\mathcal U$. First, from \eqref{rel1} we obtain
$$\theta= \frac{\pi Y}{2(X+Y)},$$
and from \eqref{rel2} we have
$$\rho=\frac ra (X+Y+Z)=\pi^{-1/3}(X+Y+Z).$$
Adding \eqref{scX} and \eqref{scY}, after some more calculations we obtain
$$ \sin \frac \varphi 2=\frac {X+Y}{\sqrt 2(X+Y+Z)},$$
and further
 $$\cos\varphi= \frac{Z(2X+2Y+Z)}{(X+Y+Z)^2},\quad \sin \varphi=\frac{X+Y}{X+Y+Z}\sqrt{2-\left(\frac{X+Y}{X+Y+Z}\right)^2}.$$
Finally, the inverse $\mathcal U^{-1}:\mathbb K \to \mathbb B^{3}$
is defined by
\begin{flalign}
x&=\pi^{-1/3}(X+Y)\sqrt{2-\left(\frac{X+Y}{X+Y+Z}\right)^2}\;\cos\frac{\pi Y}{2(X+Y)},\label{fx}\\
y&=\pi^{-1/3}(X+Y)\sqrt{2-\left(\frac{X+Y}{X+Y+Z}\right)^2}\;\sin\frac{\pi Y}{2(X+Y)},\label{fy}\\
z&=\pi^{-1/3}\frac{Z(2X+2Y+Z)}{(X+Y+Z)}.\label{fz}
\end{flalign}
for $(X,Y,Z)\in \mathbb K_0^+$, and for the other seven octants the formulas can be calculated as in \eqref{prelsim}.
\section{Uniform and refinable grids of the regular octahedron and of the ball}\label{section4}
In this section we construct a uniform refinement of  the regular
octahedron $\mathbb K$ of volume $\mbox{vol}(\mathbb K)$, more
precisely a subdivision of $\mathbb K$ into 64 cells of two
shapes, each of them having the volume $\mbox{vol}(\mathbb K)/64$.
This subdivision can be repeated for each of the 64 small cells,
the resulting $64^2$ cells of volume $\mbox{vol}(\mathbb K)/64^2$
being of one of the two types from the first refinement. Next, the
volume preserving map $\mathcal U$ will allow us the construction
of uniform and refinable grids of the 3D ball $\mathbb B^3$ by
transporting the octahedral uniform refinable 3D grids, and
further, the construction of orthonormal piecewise constant
wavelets on the 3D ball.

\subsection{Refinement of the octahedron}

The initial octahedron $\mathbb K$ consists in four congruent cells, each situated in one of the octants $I_i^+\cup I_i^-$, $i=0,1,2,3$ (see Figure \ref{fig:elem1si2}). We will say that this type of cell is $\mathbf T_0$, the index $0$ of $\mathbf T_0$ being the coarsest level of the refinement.
For simplifying the writing we denote by $\mathbb{N}_0$ the set of positive natural numbers and by  $\mathbb{N}_n=\{1,2,\ldots,n\}$, for $n\in \mathbb N_0$.
%
\subsubsection{First step of refinement}

\begin{figure}
\begin{center}
\includegraphics[width=0.45\textwidth]{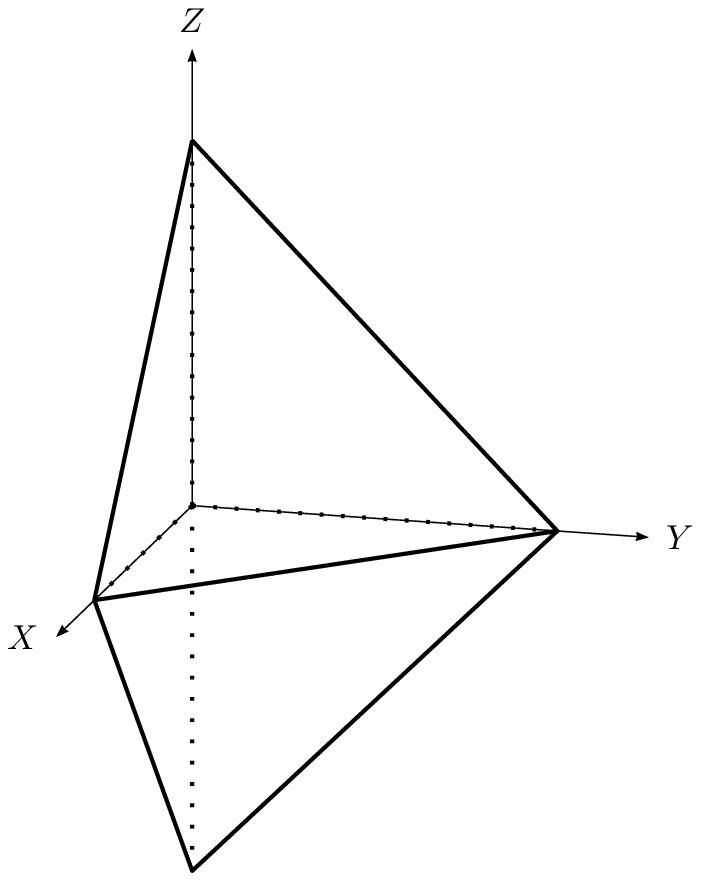}\hspace{1cm}
\includegraphics[width=0.35\textwidth]{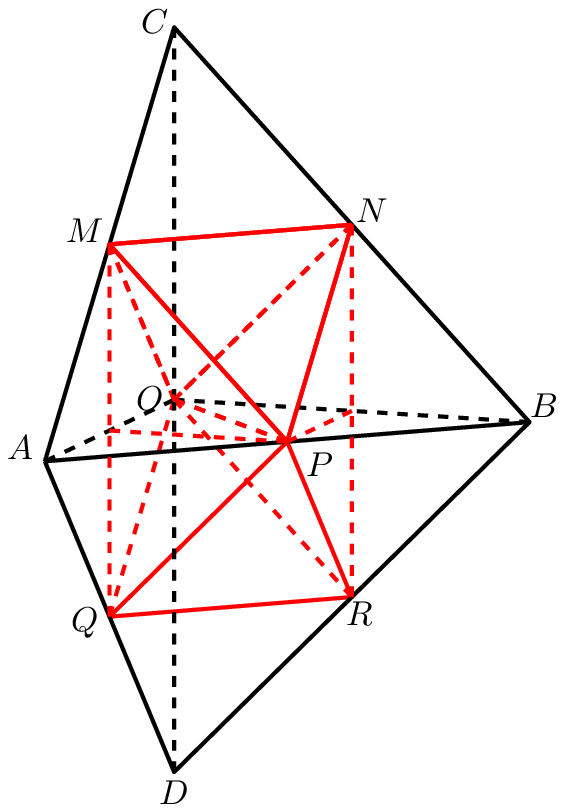}
\end{center}
\caption{Left: one of the four cells of type $\mathbf T_0$ constituting the octahedron. Right: each cell of type $\mathbf T_0$ can be subdivided into six cells of type $\mathbf T_1$ and two cells of type $\mathbf M_1$. }\label{fig:elem1si2}
\end{figure}

\begin{figure}
\begin{center}
\includegraphics*[width=0.2\textwidth]{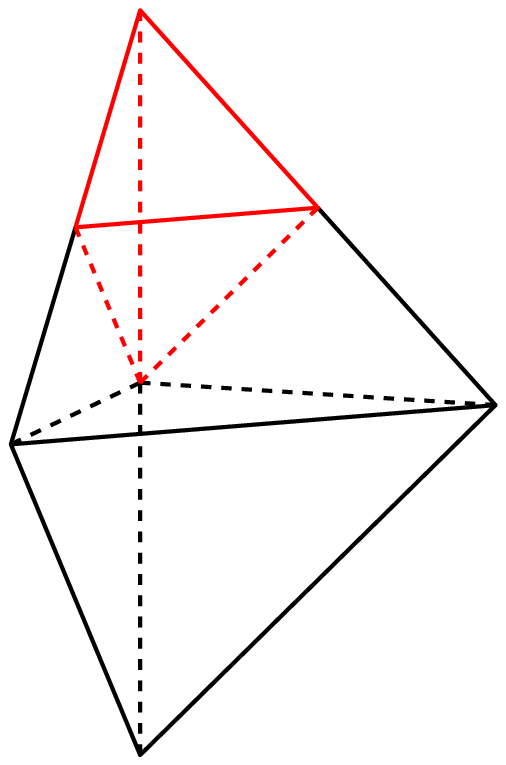}\hspace{0.01\textwidth}
\includegraphics*[width=0.2\textwidth]{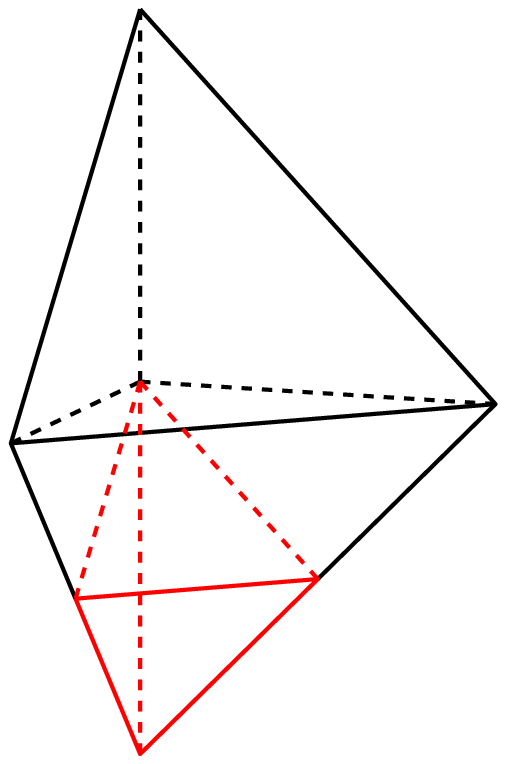}\hspace{0.01\textwidth}
\includegraphics*[width=0.2\textwidth]{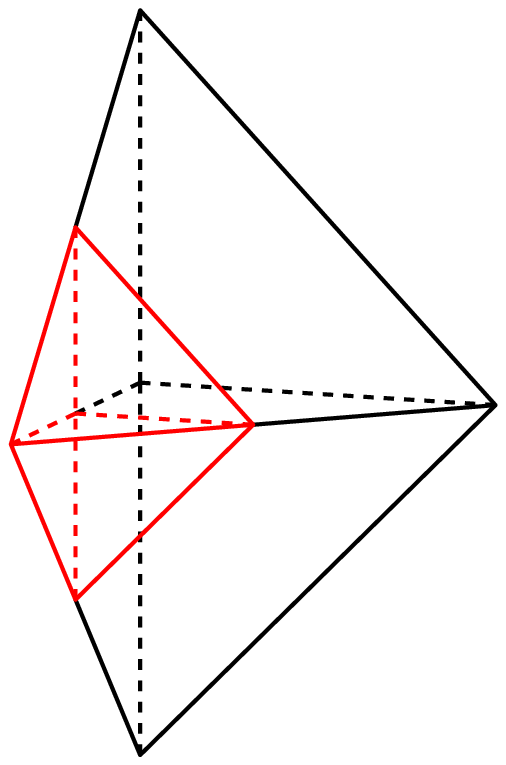}\hspace{0.01\textwidth}
\includegraphics*[width=0.2\textwidth]{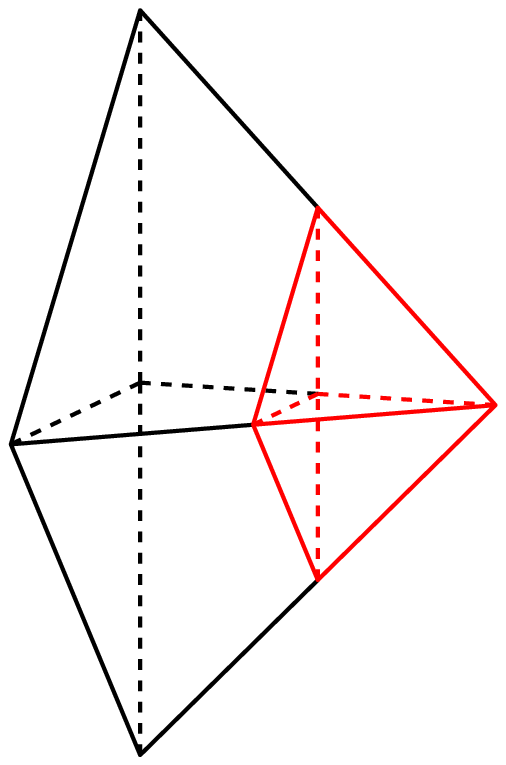}
\includegraphics*[width=0.21\textwidth]{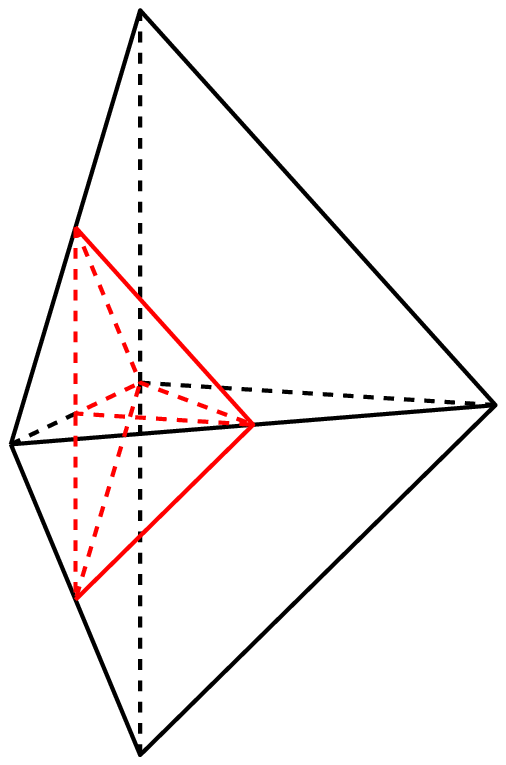}\hspace{0.01\textwidth}
\includegraphics*[width=0.21\textwidth]{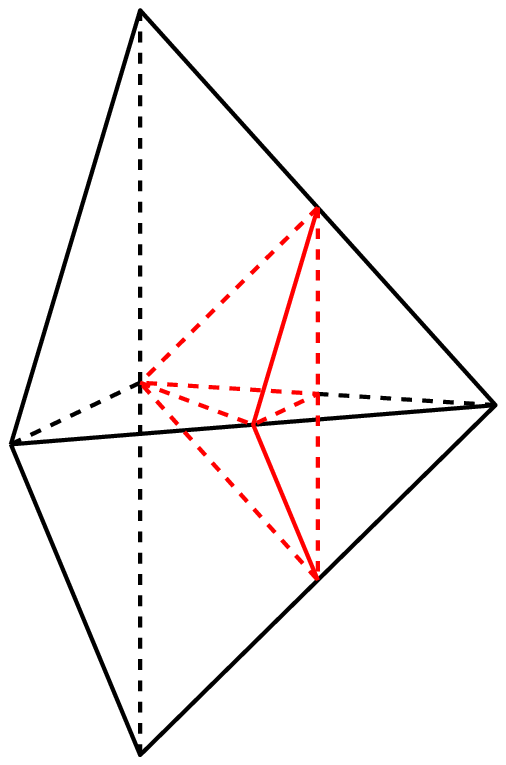}\hspace{0.01\textwidth}
\includegraphics*[width=0.21\textwidth]{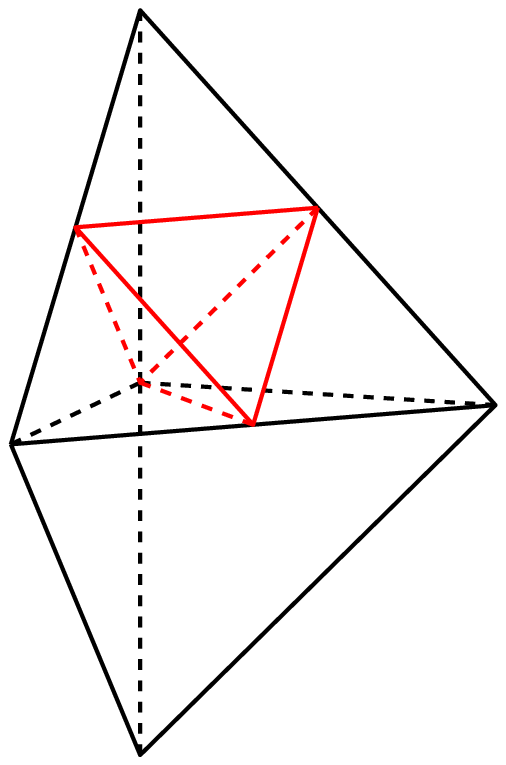}\hspace{0.01\textwidth}
\includegraphics*[width=0.21\textwidth]{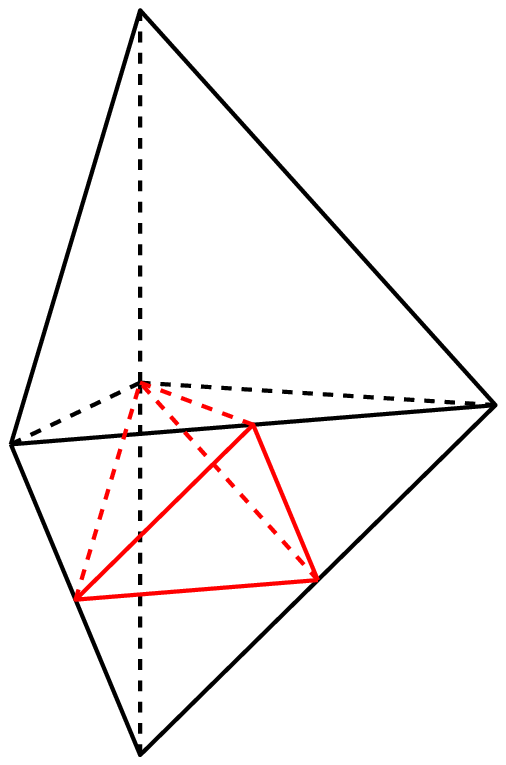}
\end{center}
\caption{The subdivision of a  $\mathbf T$ cell.}
\label{fig:desel2}
\end{figure}

\begin{figure}
\begin{center}
\includegraphics*[height=6cm]{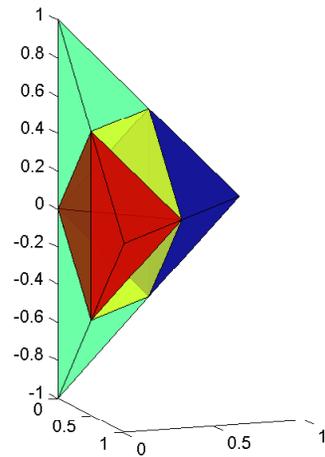}
\end{center}
\caption{The first step of the refinement: the cell $\mathbf{T}_0$
is divided into two cells of type $\mathbf{M}_1$ (yellow) and six
cells of type $\mathbf{T}_1$: two red, two blue and two
green. }\label{fig:desel2plin}
\end{figure}

The cell $\mathbf T_0=(ABCD)\in I_0^+ \cup I_0^-$, with $A(a,0,0)$, $B(0,a,0)$,
$C(0,0,a)$, $D(0,0,-a)$ (see Figure \ref{fig:elem1si2}), will be
subdivided into eight smaller cells having the same volume, as
follows: we take the mid-points $M,N,P,Q,R$ of the edges $AC$,
$BC$, $AB$, $AD$, $BD$, respectively. Thus, one obtains $t_1=6$
cells of type $\mathbf T_1$ ($MQOP$, $MQAP$, $NROP$, $NRBP$,
$ODQR$ and $COMN$),  and $m_1=2$ other cells, $OMNP$ and $OPQR$,
of another type, say $\mathbf M_1$. The cells of type $\mathbf
T_1$ have the same shape with the cells $\mathbf T_0$. Their
volumes are $$ \mbox{vol}(\mathbf T_1)= \mbox{vol}(\mathbf
M_1)=\frac{\mbox{vol}(\mathbf T_0)}{8}.$$ Figures \ref{fig:desel2}, \ref{fig:desel2plin}
also show the eight cells at the first step of refinement.

Similarly we refine the other three cells situated in $I_i^+ \cup I_i^-$, $i=1,2,3$, therefore the total number of cells after the first step of refinement is 32, more precisely 24 of type $\mathbf T_1$ and 8 of type $\mathbf M_1.$ 

\subsubsection{Second step of refinement}

A cell of type $\mathbf T_1$ will be subdivided in the same way as
a cell of type $\mathbf T_0$, i.e. into six cells of type $\mathbf
T_2$ and two cells of type $\mathbf M_2$. Their volumes will be
$$\mbox{vol}(\mathbf T_2)=\mbox{vol}(\mathbf M_2)=\frac{\mbox{vol}
(\mathbf T_0)}{8^2}.$$ Therefore, from the subdivision of the 6
cells of type $\mathbf T_1$ we have 36 cells of type $\mathbf T_2$
and 12 cells of type $\mathbf M_2$.

\begin{figure}
\begin{center}
\includegraphics*[width=0.2\textwidth]{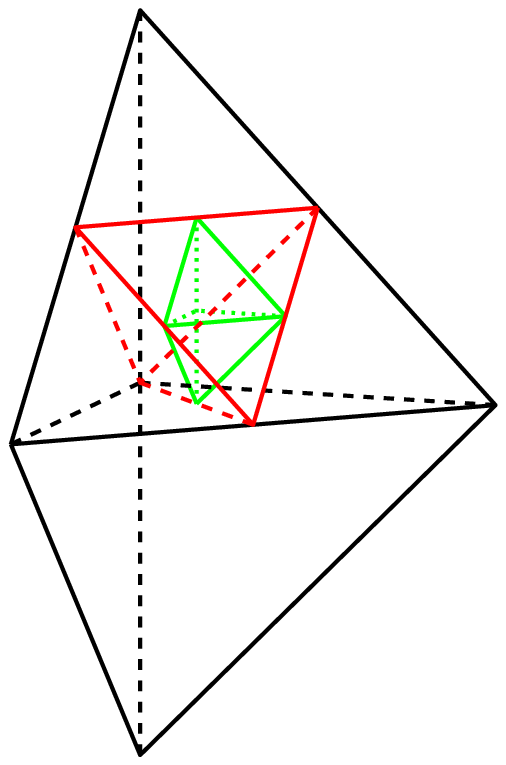}\hspace{0.01\textwidth}
\includegraphics*[width=0.2\textwidth]{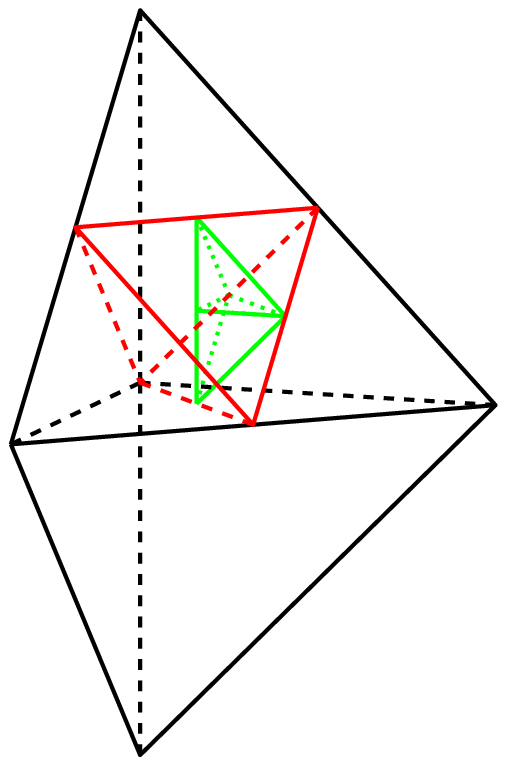}\hspace{0.01\textwidth}
\includegraphics*[width=0.2\textwidth]{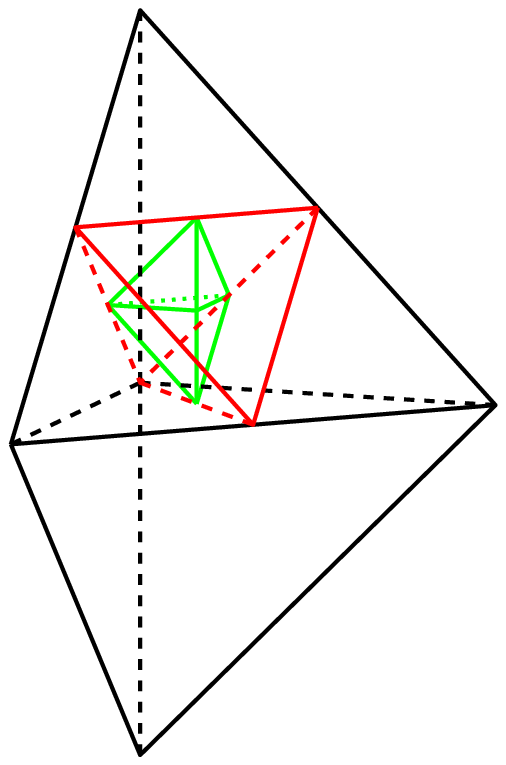}\hspace{0.01\textwidth}
\includegraphics*[width=0.2\textwidth]{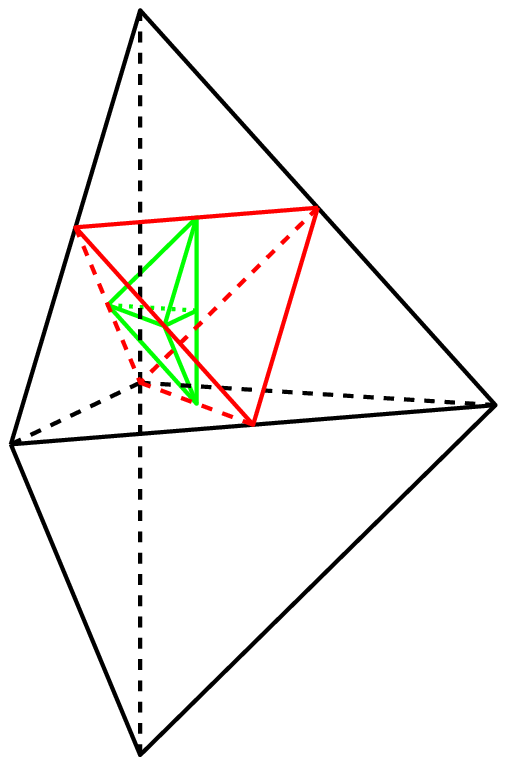}\hspace{0.01\textwidth}
\end{center}
\caption{The four cells of type $\mathbf T$ of the subdivision of a cell of type $\mathbf M$.}\label{fig:desel3}
\end{figure}

\begin{figure}
\begin{center}
\includegraphics*[width=0.4\textwidth]{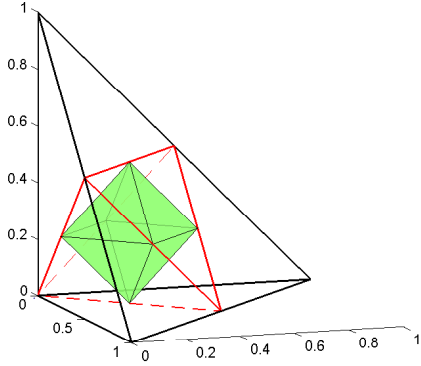}
\end{center}
\caption{The subdivision of a cell of type $\mathbf M_1$ into four cells of type $\mathbf M$: the four tetrahedrons at the corners and four cells of type $\mathbf T$ in the middle, forming an octahedron.}\label{fig:desel5}
\end{figure}

For a cell $(OMNP)$ of type $\mathbf M_1$, which is a regular tetrahedron of edge $\ell_1=a\sqrt 2/2$, we take the mid-points of the six edges (see Figures \ref{fig:desel3} and \ref{fig:desel5}). This will give four cells of type $\mathbf T_2$ in the middle and four cells of type $\mathbf M_2$, i.e regular tetrahedrons of edge $\ell_2=a\sqrt 2/2^2$. From the subdivision of the two cells of type $\mathbf M_1$ we have 8 cells of type  $\mathbf T_2$ and 8 cells of type  $\mathbf M_2$.

In conclusion, the second step of subdivision yields in $I_0^+\cup I_0^-$ $t_2= 44$
cells of type  $\mathbf T_2$ and $m_2=20$ cells of type $\mathbf
M_2$, each having the volume $\mbox{vol}(\mathbf T_0)/64,$ therefore the total number of cells after the second refinement will be $4\cdot 8^2$, more precisely 76 of type $\mathbf T_2$ and 80 of type $\mathbf M_2$.

\subsubsection{The general step of refinement}
Let $m_j$ and $t_j$ denote the numbers of cells of type $\mathbf M_j$ and $\mathbf T_j$, respectively, resulted at the step $j$ of the subdivision, starting from one cell of type $T_0$.  At this step, each of the $t_{j-1}$ cells of type  $\mathbf T_{j-1}$ is subdivided into 6 cells of type $\mathbf T_{j}$ and 2 cells of type $\mathbf M_{j}$, and each of the $m_{j-1}$ cells of type $\mathbf M_{j-1}$ is subdivided into 4 cells of type $\mathbf T_{j}$ and 4 cells of type $\mathbf M_{j}.$ This implies
\begin{eqnarray*}
t_{j}=6t_{j-1}+4m_{j-1},\\
m_{j}=2t_{j-1}+4m_{j-1},
\end{eqnarray*}
or
\begin{equation*}
\binom {t_{j}}{m_{j}}=A\binom{
          t_{j-1}}
          {m_{j-1}}
        =A^2\binom
          {t_{j-2}}
          {m_{j-2}}=\ldots= A^{j}
        \binom
          {t_{0}}
          {m_{0}},
\end{equation*}
with $t_0=1$, $m_0=0$ and $A=\left(
        \begin{array}{cc}
          6 & 4 \\
          2 & 4 \\
        \end{array}
      \right).$
After some calculations we obtain 
$$
A^j=\frac 13 \left(
      \begin{array}{cc}
         2^j (2^{2 j+1}+1)
  & 2^{j+1} (2^{2j}-1) \\
   2^j (2^{2j}-1) & 2^j (2^{2j}+2)\\
      \end{array}
    \right),\mbox{ whence}
$$
$$
t_j= \frac {2^j}3 (2^{2 j+1}+1),\quad m_j= \frac {2^j}3
(2^{2j}-1),
$$
the total number of cells of $\mathbb K_1^+\cup\mathbb K_1^-$ at
step $j$ being $t_j+m_j=8^{j},$ and $4\cdot8^j$ for the whole
octahedron $\mathbb K$. Each of the cells of type $\mathbf T_j$
and $\mathbf M_j$ has the volume $\mbox{vol}(\mathbf T_0)/8^j.$


\begin{figure}
\centering
\includegraphics*[width=0.3\textwidth]{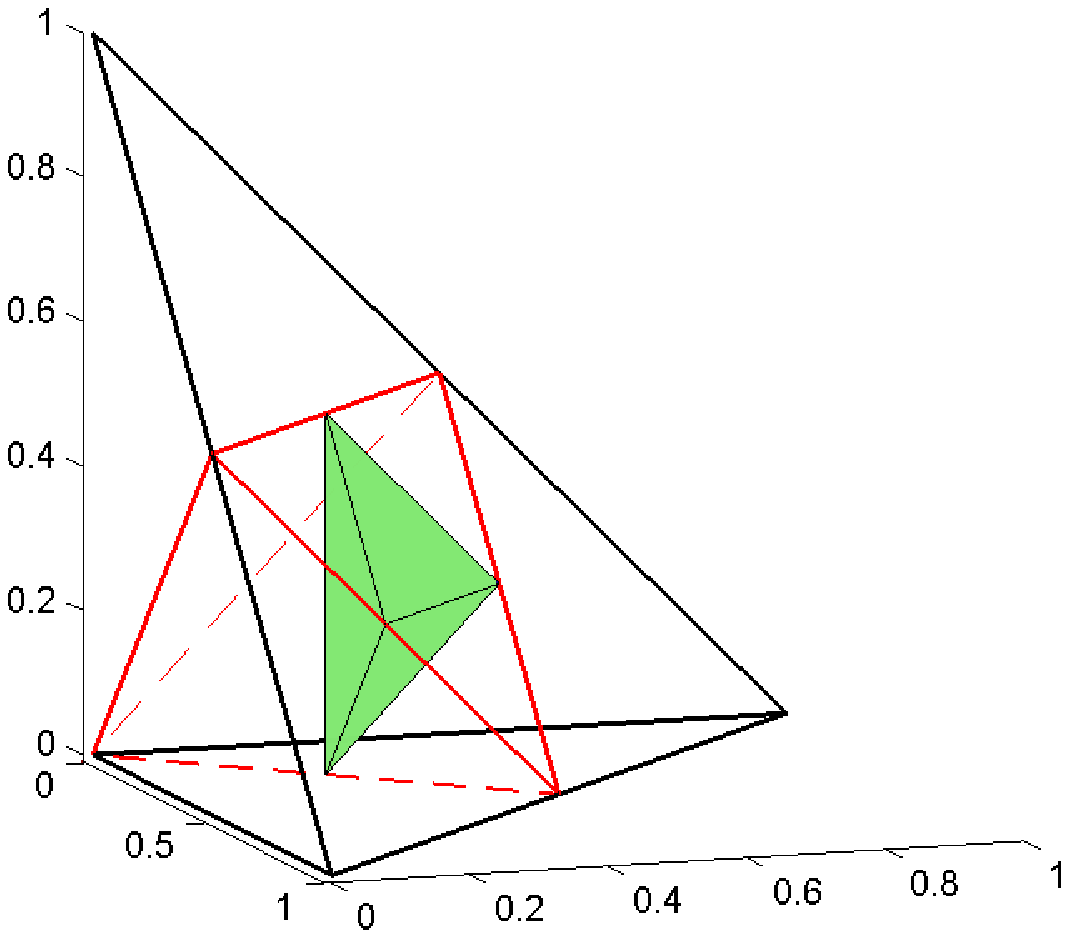}
\includegraphics*[width=0.3\textwidth]{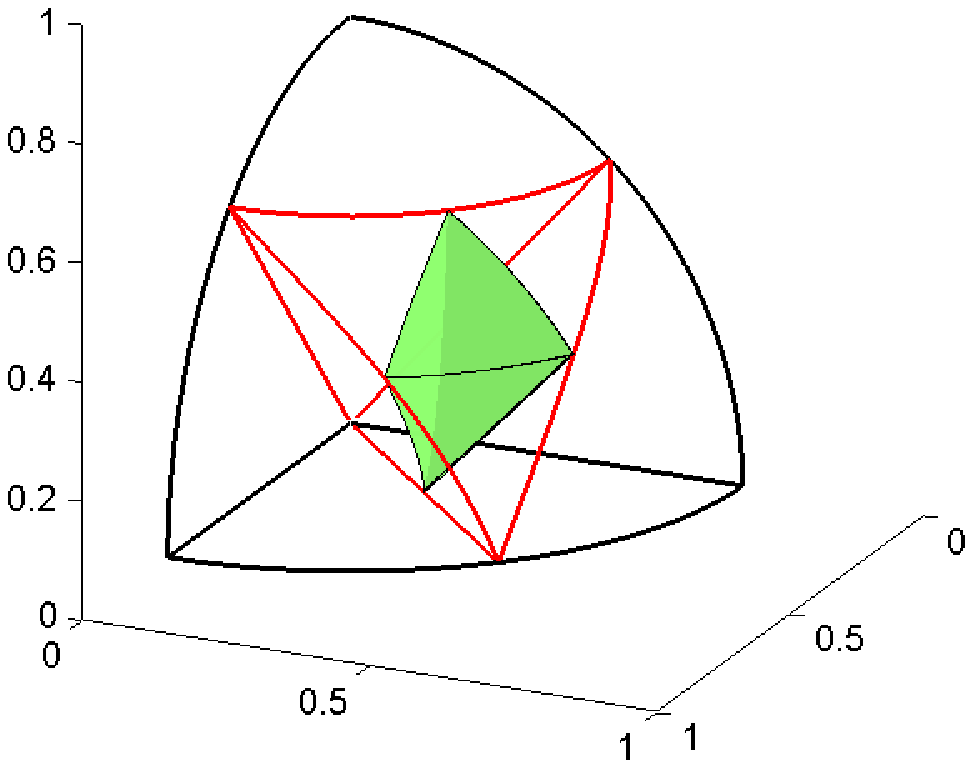}
\includegraphics*[width=0.3\textwidth]{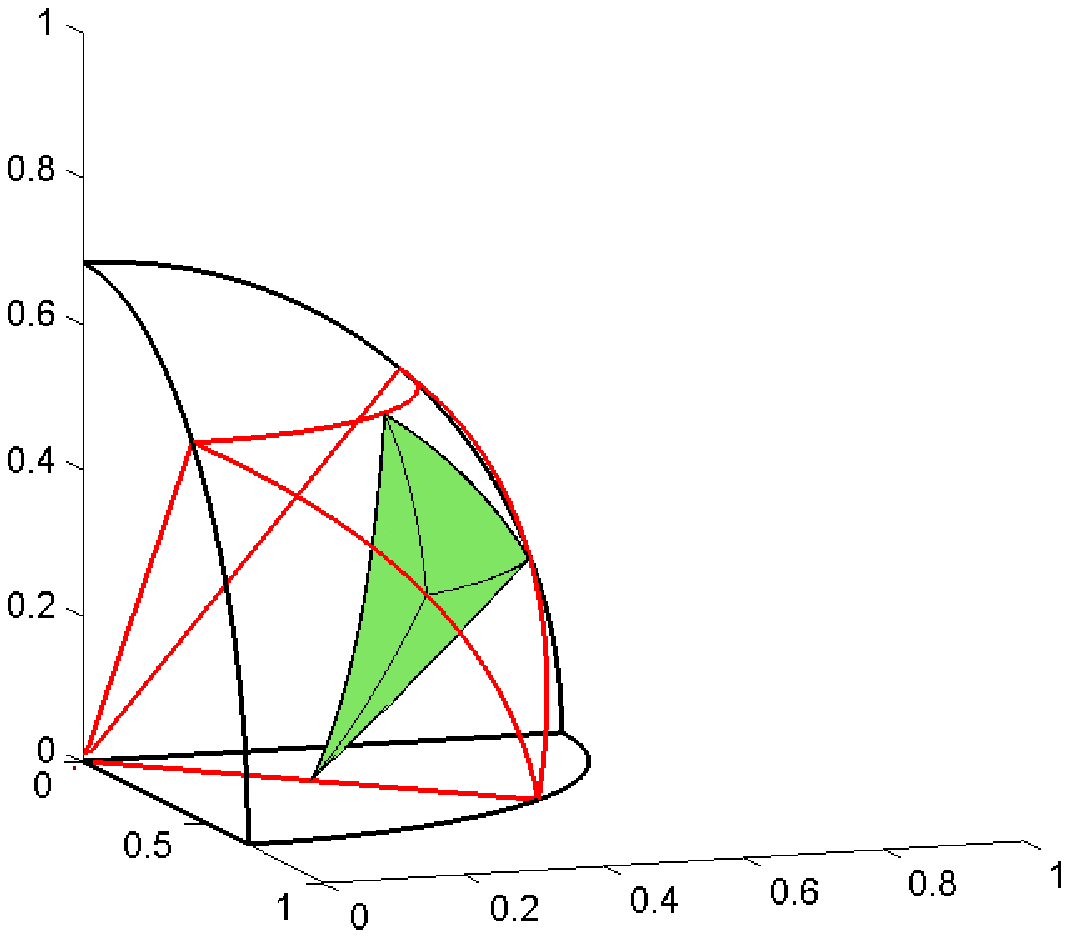}
\caption{Left: a cell of $\mathbf M$ in red and a cell of $\mathbf T$ type in green from the octahedron Middle and right: the corresponding cells of the ball.}\label{fig:desel3imag}
\end{figure}
\subsection{Implementation issues}
Every cell of the polyhedron is identified by the coordinates of its four vertices. We have two types of cells, which will be denoted  $\mathbf T$ and  $\mathbf M$.

A cell of type $\mathbf{T}$ has the same coordinates  $x$ and $y$ for the first two vertices. The $z$ coordinate of the first vertex is greater than the $z$ coordinate of the second vertex and the mean value of these $z$ coordinates gives the value of the $z$ coordinate of the third and fourth vertices of  $\mathbf T$. 

A cell of type $\mathbf M$ has two pairs of vertices at the same altitude (the same value of the $z$ coordinate).

 At every step of refinement, every cell $\mathbf T$ is divided into 6 cells of type $\mathbf T$ and two cells of type  $\mathbf M$. Suppose $[\mathbf p_1,\mathbf p_2,\mathbf p_3,\mathbf p_4]$ is the array giving the coordinates of the four vertices of a T cell. The coordinates of the vertices of the next level cells are computed as follows
\begin{flalign*}
\text{next level cell number }1:& \quad\frac{1}{2}[\mathbf p_1+\mathbf p_1,\mathbf p_1+\mathbf p_2,\mathbf p_1+\mathbf p_3,\mathbf p_1+\mathbf p_4],\\
\text{next level cell number }2:& \quad\frac{1}{2}[\mathbf p_2+\mathbf p_1,\mathbf p_2+\mathbf p_2,\mathbf p_2+\mathbf p_3,\mathbf p_2+\mathbf p_4],\\
\text{next level cell number }3:& \quad\frac{1}{2}[\mathbf p_3+\mathbf p_1,\mathbf p_3+\mathbf p_2,\mathbf p_3+\mathbf p_3,\mathbf p_3+\mathbf p_4],\\
\text{next level cell number }4:& \quad\frac{1}{2}[\mathbf p_4+\mathbf p_1,\mathbf p_4+\mathbf p_2,\mathbf p_4+\mathbf p_3,\mathbf p_4+\mathbf p_4],\\
\text{next level cell number }5:& \quad\frac{1}{2}[\mathbf p_1+\mathbf p_3,\mathbf p_2+\mathbf p_3,\mathbf p_3+\mathbf p_4,\mathbf p_1+\mathbf p_2],\\
\text{next level cell number }6:& \quad\frac{1}{2}[\mathbf p_1+\mathbf p_4,\mathbf p_2+\mathbf p_4,\mathbf p_1+\mathbf p_2,\mathbf p_3+\mathbf p_4],\\
\text{next level cell number }7:& \quad\frac{1}{2}[\mathbf p_1+\mathbf p_2,\mathbf p_1+\mathbf p_3,\mathbf p_1+\mathbf p_4,\mathbf p_3+\mathbf p_4],\\
\text{next level cell number }8:& \quad\frac{1}{2}[\mathbf p_1+\mathbf p_2,\mathbf p_2+\mathbf p_3,\mathbf p_2+\mathbf p_4,\mathbf p_3+\mathbf p_4].
\end{flalign*}
The cells 1--6 are of type $\mathbf{T}$ and the cells 7 and 8 are of type $\mathbf{M}$ (see Figure \ref{fig:desel2}).
\par Every cell  $\mathbf M$ consists in 4 cells of type $\mathbf T$ and 4 cells of type $\mathbf M$. Suppose $[\mathbf p_1,\mathbf p_2,\mathbf p_3,\mathbf p_4]$ is the array giving the coordinates of the four vertices of the cell  $\mathbf M$ and let $\mathbf p_k=\left(p_{kx}, p_{ky},p_{kz}\right)$
, $k=1,2,3,4$. We rearrange these four vertices in ascending order with respect to the $z$ coordinate. Let  $[\mathbf q_1,\mathbf q_2,\mathbf q_3,\mathbf q_4]$ be the vector $[\mathbf p_1,\mathbf p_2,\mathbf p_3,\mathbf p_4]$ sorted ascendingly with respect to the $z$ coordinate of the vertices, i.e. $q_{1z}\leq q_{2z}\leq q_{3z}\leq q_{4z}$. Similarly, let $[\mathbf r_1,\mathbf r_2,\mathbf r_3,\mathbf r_4]$ be the rearrangement of vertices $\mathbf p_1,\dots ,\mathbf p_4$ such that $r_{1x}\leq r_{2x}\leq r_{3x}\leq r_{4x}$. Let, also, $[s_1,s_2,s_3,s_4]$ be the array of rearranged vertices with respect to the $y$ coordinate in ascending order. The coordinates of the vertices of the cells at the next level are computed as follows:
\begin{flalign*}
\text{next level cell number }1:& \quad\frac{1}{2}[\mathbf q_3+\mathbf q_4,\mathbf q_1+\mathbf q_2,\mathbf r_3+\mathbf r_4,\mathbf s_3+\mathbf s_4]\\
\text{next level cell number }2:& \quad\frac{1}{2}[\mathbf q_3+\mathbf q_4,\mathbf q_1+\mathbf q_2,\mathbf s_3+\mathbf s_4,\mathbf r_1+\mathbf r_2]\\
\text{next level cell number }3:& \quad\frac{1}{2}[\mathbf q_3+\mathbf q_4,\mathbf q_1+\mathbf q_2,\mathbf r_1+\mathbf r_2,\mathbf s_1+\mathbf s_2]\\
\text{next level cell number }4:& \quad\frac{1}{2}[\mathbf q_3+\mathbf q_4,\mathbf q_1+\mathbf q_2,\mathbf s_1+\mathbf s_2,\mathbf r_3+\mathbf r_4]\\
\text{next level cell number }5:& \quad\frac{1}{2}[\mathbf p_1+\mathbf p_1,\mathbf p_1+\mathbf p_2,\mathbf p_1+\mathbf p_3,\mathbf p_1+\mathbf p_4]\\
\text{next level cell number }6:& \quad\frac{1}{2}[\mathbf p_2+\mathbf p_1,\mathbf p_2+\mathbf p_2,\mathbf p_2+\mathbf p_3,\mathbf p_2+\mathbf p_4]\\
\text{next level cell number }7:& \quad\frac{1}{2}[\mathbf p_3+\mathbf p_1,\mathbf p_3+\mathbf p_2,\mathbf p_3+\mathbf p_3,\mathbf p_3+\mathbf p_4]\\
\text{next level cell number }8:& \quad\frac{1}{2}[\mathbf p_4+\mathbf p_1,\mathbf p_4+\mathbf p_2,\mathbf p_4+\mathbf p_3,\mathbf p_4+\mathbf p_4].
\end{flalign*}
To verify whether a point $\mathbf p=(p_x,p_y,p_z)$  is inside a cell with vertices $\left[\mathbf p_1, \mathbf p_2,\mathbf p_3,\mathbf p_4 \right]$, we compute the following numbers:
\begin{flalign*}
d_1&=\sgn\begin{vmatrix}p_{1x} & p_{2x} & p_{3x} & p_x\\p_{1y} & p_{2y} & p_{3y} & p_y\\ p_{1z} & p_{2z} & p_{3z} & p_z\\1 & 1 & 1 & 1 \end{vmatrix}, d_2=\sgn\begin{vmatrix}p_{1x} & p_{2x} & p_x & p_{4x}\\p_{1y} & p_{2y} & p_y & p_{4y}\\ p_{1z} & p_{2z} & p_z & p_{4z}\\1 & 1 & 1 & 1 \end{vmatrix}, d_3=\sgn\begin{vmatrix}p_{1x} & p_{x} & p_{3x} & p_{4x}\\p_{1y} & p_{y} & p_{3y} & p_{4y}\\ p_{1z} & p_{z} & p_{3z} & p_{4z}\\1 & 1 & 1 & 1 \end{vmatrix},\\ d_4 &=\sgn\begin{vmatrix}p_{x} & p_{2x} & p_{3x} & p_{4x}\\p_{y} & p_{2y} & p_{3y} & p_{4y}\\ p_{z} & p_{2z} & p_{3z} & p_{4z}\\1 & 1 & 1 & 1 \end{vmatrix}, d_5=\sgn\begin{vmatrix}p_{1x} & p_{2x} & p_{3x} & p_{4x}\\p_{1y} & p_{2y} & p_{3y} & p_{4y}\\ p_{1z} & p_{2z} & p_{3z} & p_{4z}\\1 & 1 & 1 & 1 \end{vmatrix}.
\end{flalign*}
We calculate $v=|d_1|+|d_2|+|d_2|+|d_3|+|d_4|+|d_5|$. If $|d_1+d_2+d_3+d_4+d_5|=v$, then for $v=5$ the point $\mathbf p$ is in the interior of the cell, for $v=4$ the point $ \mathbf p$ is on one of the faces of the cell, for $v=3$ the point $\mathbf p$ is situated on one of the edges of the cell, and for $v=2$ the point $ \mathbf p$ is one of the vertices of the cell. If $|d_1+d_2+d_3+d_4+d_5|\neq v$, the point $\mathbf p$ is located outside the cell. Since the vertices $\mathbf p_k$ are different we have $v\geq 2$.
\begin{figure}
\begin{center}
\includegraphics*[width=0.4\textwidth]{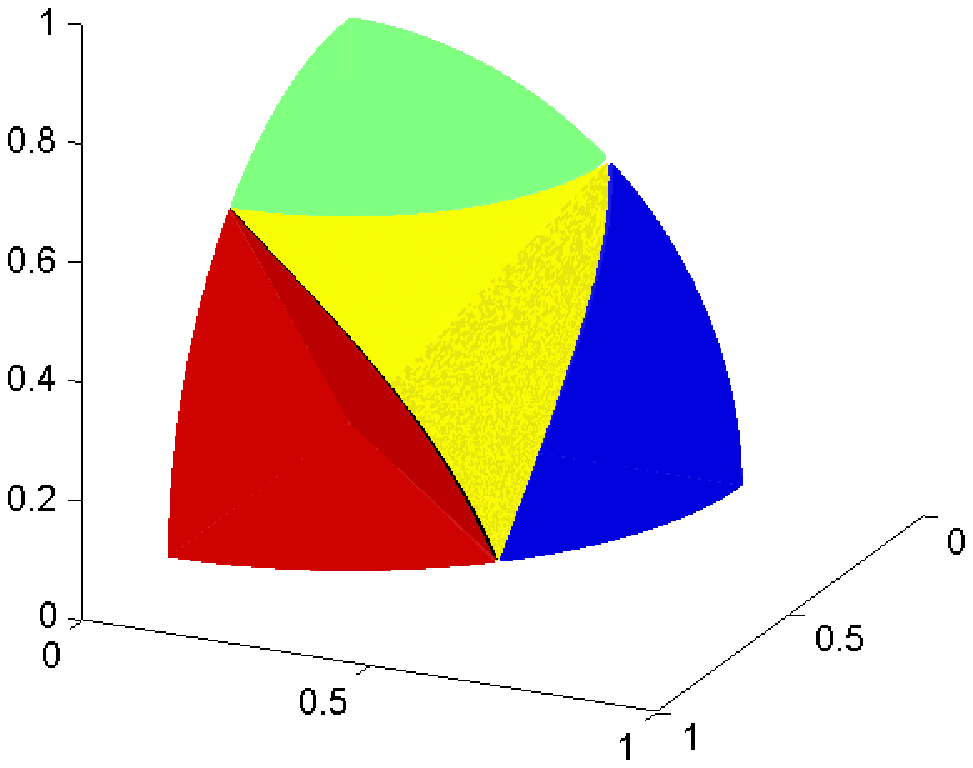}\hspace{1cm}
\includegraphics*[width=0.4\textwidth]{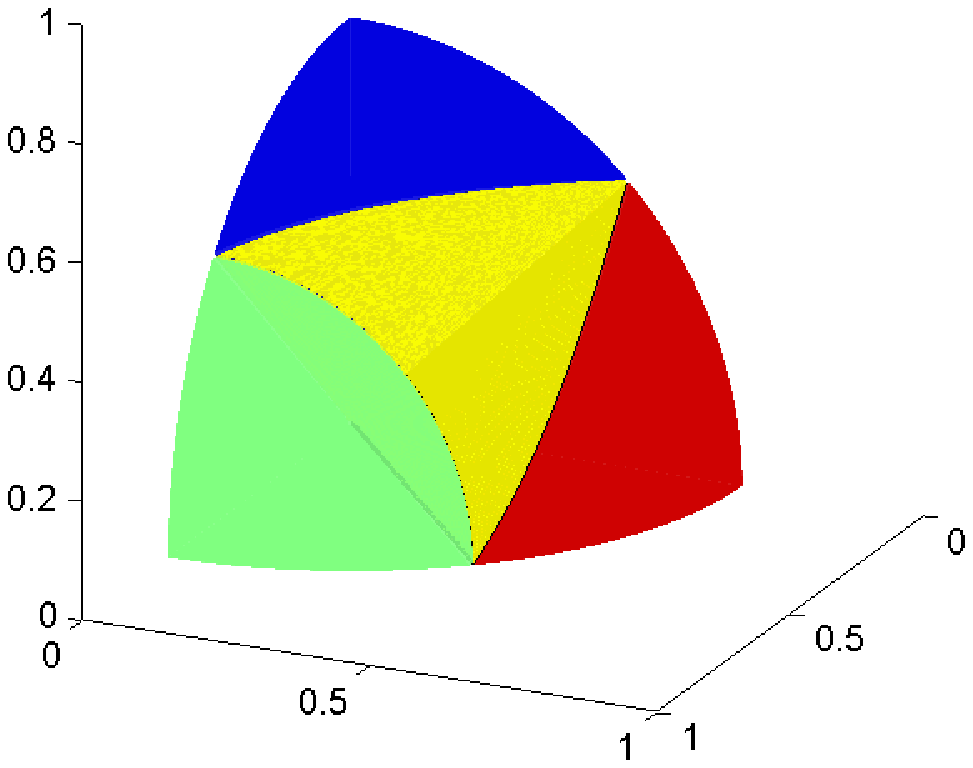}\hspace{1cm}
\end{center}
\caption{Left: the image on the ball of the positive octant; Right: the same image rotated.}\label{fig:desel2imagpoz}
\end{figure}

\begin{figure}
\centering
\includegraphics*[width=0.4\textwidth]{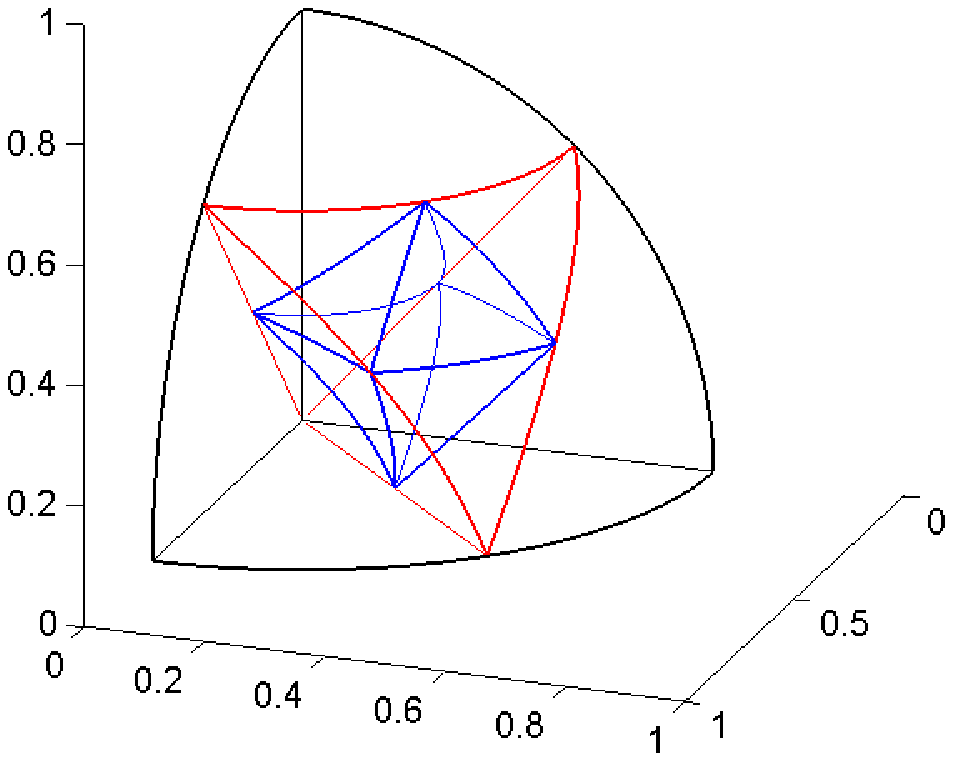}\hspace{1cm}
\includegraphics*[width=0.4\textwidth]{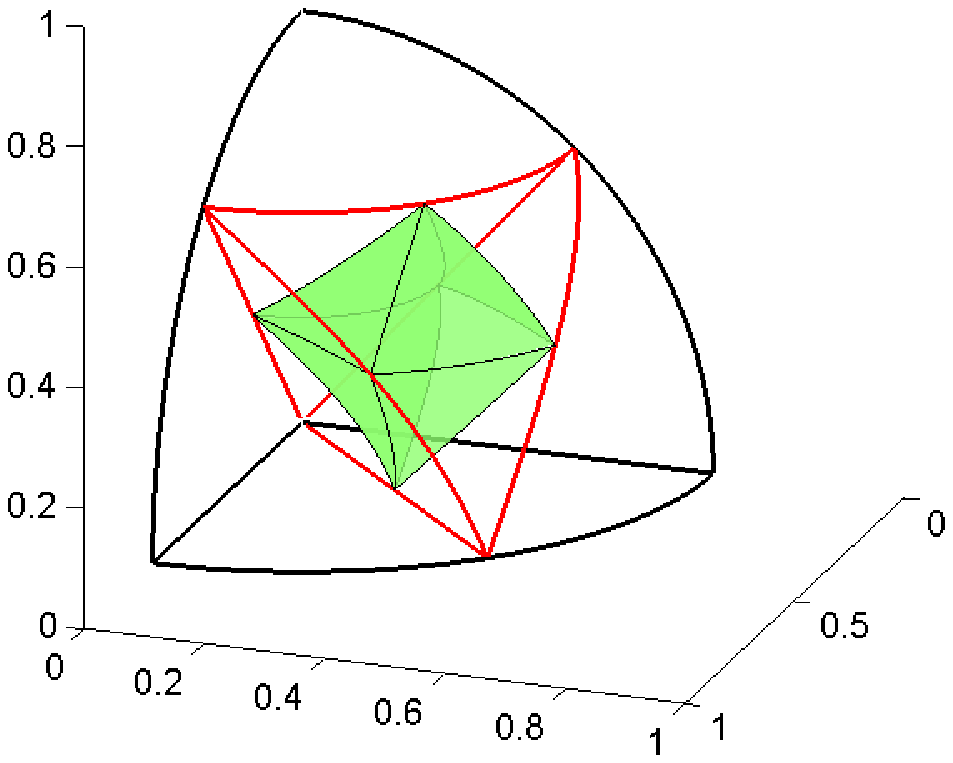}
\caption{The image on the ball of the cells of the octahedron corresponding to Figure \ref{fig:desel5}.}
\label{fig10}
\end{figure}

\subsection{Uniform and refinable grids of the ball $\mathbb B^{3}$}
If we transport the uniform and refinable grid on $\mathbb K$ onto the ball $\mathbb B^3$ using the volume preserving map $\mathcal U^{-1}$, we obtain a uniform and refinable grid of $\mathbb B^3$.  Figures \ref{fig:desel3imag},  \ref{fig:desel2imagpoz} and  \ref{fig10} show the images on $\mathbb B^3$  of different cells of $\mathbb K$.

Besides the multiresolution analysis and wavelet bases, which will be constructed in Section \ref{section5}, another useful application is the construction of a uniform sampling of the rotation group $SO(3)$, by calculations similar to the ones in \cite{marc}. This will be subject of a future paper.

\section{Multiresolution analysis and piecewise constant orthonormal wavelet bases of $L^2(\mathbb K)$ and $L^2(\mathbb B^3)$}\label{section5}
Let $\mathcal D=\mathcal D^0=\{D_1,D_2,D_3,D_4\}$ be the
decomposition of the domain $\mathbb K$ considered in Section 4.1,
consisting in four congruent domains (cells) of type $\mathbf
T_0$. For $D\in \mathcal D,$ let $\mathcal R_D$ denote the set of
the eight refined domains, constructed in Section 4.1.1. The set
$\mathcal D^1=\cup_{D\in\mathcal D^0}\mathcal R_D$ is a refinement
of $\mathcal D^0$, consisting in $4\cdot 8$ congruent cells.
Continuing the refinement process as we described in Section 4, we
obtain a decomposition $\mathcal D^j$ of $\mathbb K,$ for $j\in
\mathbb N_0,$ $|\mathcal D^j|=4\cdot 8^j$.

For a fixed $j\in \mathbb N_0$ we assign to each domain $D_k^j\in
\mathcal D^j$, $k\in \mathcal N_j:=\mathbb N_{4\cdot 8^j}$, the
function $\varphi_{D_k^j}:\mathbb K \to \mathbb R,$
$$
\varphi_{D_k^j}=(2\sqrt 2)^j\frac 2 {\sqrt{\mbox{vol} (\mathbb
K)}}\chi_{D_k^j} ,
$$
where $\chi_{D_k^j}$ is the characteristic function of the domain
$D_k^j$. Then we define the spaces of functions
$V^j=\mbox{span}\{\varphi_{D_k^j},\ k\in \mathcal N_j\} $ of
dimension $4\cdot 8^j$, consisting of piecewise constant functions
on the domains of $\mathcal D^j.$ Moreover, we have
$\|\varphi_{D_k^j}\|_{L^2(\mathbb K)}=1$, the norm being the usual $2$-norm of the
space $L^2(\mathbb K)$. For $A^j\in \mathcal D^j=\{D_k^j,\ j\in \mathcal N_j\}$, let
$A_k^{j+1},$ $k\in \mathbb N_8$, be the refined subdomains obtained
from $A^j$. One has
$$
\varphi_{A^j}=\frac 1{2\sqrt 2}\left( \varphi_{A_1^{j+1}}+
\varphi_{A_2^{j+1}}+\ldots+ \varphi_{A_8^{j+1}}\right),
$$
in $L^2(\mathbb K)$, equality which implies  the inclusion $V^j
\subseteq V^{j+1}$, for all $j\in \mathbb N_0$. With respect to the
usual inner product $\langle \cdot,\cdot \rangle_{L^2(\mathbb K)}$, the spaces
$V^j$ are Hilbert spaces, with the corresponding usual 2-norm $\| \cdot \|_{L^2(\mathbb K)}$.
In conclusion, the sequence of subspaces $V^j$ has the following properties:
\begin{enumerate}
\item $V^j\subseteq V^{j+1} \mbox{ for all }j\in \mathbb N_0$,
\item $\mbox{clos}_{L^2(\mathbb K)} \bigcup_{j=0}^\infty V^j=L^2(\mathbb K)$,
\item The set $\{\varphi_{D_k^j},\ k\in \mathcal N_j \}$ is an orthonormal basis of the space $V^j$ for each $j\in \mathbb N_0$,
\end{enumerate}
i.e. the sequence $\{V^j,\ j\in \mathbb N_0 \}$ constitutes a \emph{multiresolution analysis} of the space $L^2(\mathbb K)$.
 Let
$W^j$ denote the orthogonal complement of the coarse space $V^j$
in the fine space $V^{j+1},$ so that
$$
V^{j+1}=V^j \oplus W^j.
$$
The dimension of $W^j$ is $\mbox{dim }W^j=28 \cdot 8^j$. The
spaces $W^j$ are called \textit{wavelet spaces} and their elements are called \emph{wavelets}. In the following
we construct an orthonormal basis of $W^j$. To each domain $A^j\in \mathcal D^j$,
seven wavelets supported on $D^j$ will be associated in the
following way:
$$
\psi_{A^j}^{\ell}=a_{\ell 1}\varphi_{A_1^{j+1}}+a_{\ell 2}\varphi_{A_2^{j+1}}+\ldots+a_{\ell 8}\varphi_{A_8^{j+1}},
\ \mbox{for }\ell\in \mathbb N_7,
$$
with $a_{\ell j}\in \mathbb R,$ $\ell \in \mathbb N_7,$ $j\in \mathbb N_8.$ We
have to find conditions on the coefficients $a_{\ell j}$ which ensure
that the set $\{ \psi_{A^j}^{\ell},\ \ell\in \mathbb N_7,\ A^j\in \mathcal
D^j\}$ is an orthonormal basis of $W^j$. First we must have
\begin{equation}\label{condort}
\langle \psi_{A^j}^{\ell}, \varphi_{S^j}\rangle=0, \mbox{ for
}\ell\in \mathbb N_7 \mbox{ and }A^j,S^j\in \mathcal D^j.
\end{equation}
If $A^j \neq S^j$, the equality is immediate, since $\mbox{supp
}\psi_{A^j}^{\ell}\subseteq \mbox{supp }\varphi_{A^j}$ and $\mbox{supp
}\varphi_{A^j} \cap \mbox{supp }\varphi_{S^j}$ is either empty or
an edge, whose measure is zero. If $A^j=S^j$, evaluating the inner
product \eqref{condort} we obtain
\begin{eqnarray*}
\langle \psi_{A^j}^\ell, \varphi_{S^j}\rangle&=& \langle
a_{\ell 1}\varphi_{A_1^{j+1}}+a_{\ell 2}\varphi_{A_2^{j+1}}+\ldots+a_{\ell 8}\varphi_{A_8^{j+1}},\varphi_{A^j}\rangle
\\&=&\frac
1{2\sqrt 2}(a_{\ell 1}+a_{\ell 2}+\ldots+a_{\ell 8}).
\end{eqnarray*}
Then, each of the orthogonality conditions
$$
\langle \psi_{A^j}^\ell, \psi_{A^j}^{\ell \,'} \rangle=\delta_{\ell \ell\, '},
\mbox{ for all }A^j\in \mathcal D^j,
$$
is equivalent to
$a_{\ell \,'1}a_{\ell1}+a_{\ell\, '2}a_{\ell2}+\ldots+a_{\ell \,'8}a_{\ell8}=\delta_{\ell\ell\,'},$
$\ell,\ell\,'\in \mathbb N_7$. In fact, one requires the orthogonality of
the $8\times 8$ matrix $M=\left(a_{ij} \right)_{i,j}$ with the
entries of the first row equal to $1/(2\sqrt 2)$.

A particular case was considered in \cite{Argentina}, where the
authors divide a tetrahedron into eight small tetrahedrons of the
same area using Bey's method and for the construction of the
orthonormal wavelet basis they take the Haar matrix
$$
\frac 1{2\sqrt 2}\left(%
\begin{array}{cccccccc}
  1 & 1 & 1 & 1 & 1 & 1 & 1 & 1 \\
   1 & 1 & 1 & 1 & -1 & -1 & -1 & -1 \\
  1 & 1 & -1 & -1 & 0 & 0 & 0 & 0 \\
  0 & 0 & 0 & 0 & 1 & 1 & -1 & -1 \\
  1 & -1 & 0 & 0 & 0 & 0 & 0 & 0 \\
  0 & 0 & 1 & -1 & 0 & 0 & 0 & 0 \\
  0 & 0 & 0 & 0 & 1 & -1 & 0 & 0 \\
  0 & 0 & 0 & 0 & 0 & 0 & 1 & -1 \\
\end{array}%
\right)
$$
Alternatively, we can consider the symmetric orthogonal matrix
$$\left(%
\begin{array}{cccccccc}
  c & c & c & c & c & c & c & c \\
  c & a & b & b & b & b & b & b \\
  c & b & a & b & b & b & b & b \\
  c & b & b & a & b & b & b & b \\
  c & b & b & b & a & b & b & b \\
  c & b & b & b & b & a & b & b \\
  c & b & b & b & b & b & a & b \\
  c & b & b & b & b & b & b & a \\
\end{array}
\right),$$ with $$a= \frac{\pm 24-\sqrt 2}{28},\ b=\frac {\mp
4-\sqrt 2}{28},\ c=\frac 1{2\sqrt 2},$$ or the tensor product $H\otimes H \otimes H$ of the matrix
$$
H=\frac 1{\sqrt 2}\left(
                  \begin{array}{cc}
                    1 & 1\\
                    1 & -1 \\
                  \end{array}
                \right), \mbox {which is}
$$
\begin{eqnarray*}
\frac 1{2\sqrt 2}\left(%
\begin{array}{cccccccc}
  1 & 1 & 1 & 1 & 1 & 1 & 1 & 1 \\
   1 & -1 & 1 & -1 & 1 & -1 & 1 & -1 \\
  1 & 1 & -1 & -1 & 1 & 1 & -1 & -1 \\
  1 &-1 & -1 & 1 & 1 & -1 & -1 & 1 \\
  1 & 1 & 1 & 1 & -1 & -1 & -1 & -1 \\
  1 & -1 & 1 & -1 & -1 & 1 & -1 & 1 \\
  1 & 1 & -1 & -1 & -1 & -1 & 1 & 1 \\
  1 & -1 & -1 & 1 & -1 & 1 & 1 & -1 \\
\end{array}%
\right)
\end{eqnarray*}
 or, more
general, we can generate \emph{all} orthogonal $8\times 8$
matrices with the entries of the first row equal to $1/(2\sqrt 2)$
using the method described in \cite{poprosca}, where we start with
the well known Euler's formula for the general form of a $3\times
3$ rotation matrix. It is also possible to use different orthogonal matrices
for the wavelets associated to the decomposition of the cells of
type $\mathbf T$ and $\mathbf M$.

Next, following the ideas in \cite{ACHA} we show how one can transport the above multiresolution analysis and wavelet bases on the 3D ball $\mathbb B^3$, using the volume preserving map $\mathcal U:\mathbb B^3 \to \mathbb K$ constructed in Section \ref{section3}.

We that  consider the ball $\mathbb B^3$ is given by the parametric equations
$$
\xi=\xi(X,Y,Z)=\mathcal U^{-1}(X,Y,Z)=\left(x(X,Y,Z),y(X,Y,Z),z(X,Y,Z) \right),
$$
with $(X,Y,Z)\in \mathbb K$. Since $\mathcal U$ and its inverse preserve the volume, the volume element $d\omega (\xi)$ of $\mathbb B^3$ equals the volume element $dX\,dY\,dZ=d\mathbf x$ of $\mathbb K$ (and $\mathbb R^3$). Therefore, for all $\widetilde{f},\widetilde{g}\in L^2(\mathbb B^3)$ we have
\begin{eqnarray*}
\langle \widetilde{f},\widetilde{g} \rangle _{L^2(\mathbb B^3)}&=&\int_{\mathbb B^3}\overline{\widetilde{f}(\xi)}\widetilde{g}(\xi)\, d\omega(\xi)\\
&=& \int_{\mathcal U(\mathbb B^3)}\overline{\widetilde{f}(\mathcal U^{-1}(X,Y,Z))}\, \widetilde{g}(\mathcal U^{-1}(X,Y,Z))\,dX\,dY\,dZ\\
&=& \langle \widetilde{f}\circ \mathcal U^{-1},\widetilde{g} \circ \mathcal U^{-1}\rangle _{L^2(\mathbb K),}
\end{eqnarray*}
and similarly, for all $f,g \in L^2(\mathbb K)$ we have
\begin{equation}\label{innpr}
\langle f,g \rangle_{L^2(\mathbb K)}=\langle f \circ \mathcal U,g\circ \mathcal U \rangle_{L^2(\mathbb B^3)}.
\end{equation}

If we consider the map $\Pi: L^2(\mathbb B^3)\to L^2(\mathbb K)$ induced by $\mathcal U$,
 defined by
 $$
 (\Pi \widetilde{f})(X,Y,Z)=\widetilde{f}\left( \mathcal U^{-1}(X,Y,Z) \right ), \mbox{ for all } \widetilde{f}\in L^2(\mathbb B^3),
 $$
and its  inverse $\Pi^{-1}:L^2(\mathbb K) \to L^2(\mathbb B^3),$
$$
(\Pi^{-1}f)(\xi)=f(\mathcal U(\xi)), \mbox{ for all } f\in L^2(\mathbb K),
$$
then $\Pi$ is a unitary map, that is
\begin{eqnarray}
\langle \Pi \widetilde{f},\Pi \widetilde{g} \rangle_{L^2(\mathbb K)}=\langle \widetilde{f},\widetilde{g} \rangle_{L^2(\mathbb B^3)},\label{um1}\\
\langle \Pi^{-1} {f},\Pi^{-1} {g} \rangle_{L^2(\mathbb B^3)}=\langle {f},{g} \rangle_{L^2(\mathbb K)}.\label{um2}
\end{eqnarray}

Equality \eqref{innpr} suggests us the construction of orthonormal scaling functions and wavelets defined on $\mathbb B^3$. The scaling functions $\widetilde{\varphi_{D^j_k}}:\mathbb B^3 \to \mathbb R$ will be
\begin{equation}
\widetilde{\varphi_{D^j_k}}=\varphi_{D^j_k} \circ \mathcal U=\left \{
\begin{array}{cl}
  1, & \mbox{on } \mathcal U^{-1}(D_k^j), \\
  0 ,& \mbox{in rest}.
\end{array}
\right .
\end{equation}
and the wavelets will be defined similarly,
$$
\widetilde{\psi_{A^j}^\ell}=\psi_{A^j}^{\ell} \circ \mathcal U.
$$

From equality \eqref{innpr} we can conclude that the spaces
$$
\widetilde{V^j}:=\mbox{span }\{ \widetilde{\varphi_{D_k^j}},\, k\in \mathcal N_j\}
$$
constitute a multiresolution analysis of $L^2(\mathbb B^3)$, each of the set
 $\{ \widetilde{\varphi_{D_k^j}},\, k\in \mathcal N_j\}$ being an orthonormal basis for the space $\widetilde{V^j}$. Moreover,
 the set
 $$\{ \widetilde{\psi^{\ell}_{A^j}},\, \ell \in \mathbb N_7,\, A^j\in \mathcal D_j\}$$ is an orthonormal basis of $\widetilde{W^j}$.

\end{document}